\documentclass[a4paper,11pt]{article}
\usepackage{amsmath,amscd,amsfonts,amssymb,latexsym}

\newtheorem{theorem}{Theorem}

\newtheorem{example}[theorem]{Example}

\def\Hom{{\rm Hom}}
\def\C{\mathbb{C}}
\def\Q{\mathbb{Q}}
\def\T{\mathbb{T}}

\def\O{{\cal O}}

\def\P{\mathbb{P}}

\def\dim{{\rm dim}}

\author{Andrzej Weber \break University of Warsaw}

\begin{document}
\title{\bf Computing equivariant characteristic classes of singular varieties
}
\author{
Andrzej Weber\\
\small Department of Mathematics of Warsaw University\\
\small Banacha 2, 02-097 Warszawa, Poland\\
\small aweber@mimuw.edu.pl\\
}
\date{June 2013}

\maketitle

\begin{abstract} Starting from the classical theory we describe Hirzebruch class and the related Todd genus of a complex singular  algebraic varieties. When the variety is equipped with an action of an algebraic torus we localize the Hirzebruch class at the fixed points of the action. We give some examples of computations.

 The paper is  based on the talk given in Kyoto, November 2012.
\end{abstract}

\section{Riemann-Roch theorem}

More than 150 years ego Bernhard Riemann  \cite{Ri}  proved
certain inequality, which in the contemporary language can be
stated as follows: Let $C$ be a smooth complete complex curve
(i.e. a Riemann surface) of genus $g$ and let $L$ be a line bundle
over $C$. Then
$$h^0(C;L)\geq \deg(L)+1-g\,,$$
where $h^0(C;L)$ denotes the dimension of the space of the
holomorphic sections $H^0(C;L)$ and $\deg(L)$ is the degree of the
line bundle which is equal to the integral of its first Chern class
$$\deg(L)=\int_Cc_1(L)\,.$$
Few years later a student of Riemann, Gustav Roch computed the
error term in the inequality. Now due to Serre duality we can
write
\begin{theorem}$$h^0(C;L)-h^1(C;L)=\deg(L)+1-g\,.$$\end{theorem}
This number is equal to the Euler characteristic $\chi(C;L)$ of
the line bundle $L$. The Riemann-Roch theorem  became an
indispensable tool for algebraic geometers. It was generalized to
higher dimensions by Hirzebruch: Let $X$ be a projective smooth
algebraic variety,   $E$ an algebraic vector bundle over $X$, then
\begin{theorem}
$$\chi(X;E)=\int_{X} ch(E)\cup td(TX).$$\end{theorem}
In this formula  $ch(-)$ and $td(-)$ are  expressions in Chern
classes of a vector bundle. The first ingredient of the formula is
the Chern character of a vector bundle. It satisfies two
conditions:
 \begin{itemize}
 \item for a line bundle $$ch(L)=exp(c_1(L))=1+c_1(L)+\frac{c_1(L)^2}{2!}+\frac{c_1(L)^3}{3!}+\frac{c_1(L)^4}{4!}+\dots\,,$$
 \item Chern character is additive
 $$ch(E_1\oplus E_2)=ch(E_1)+ch(E_2)\,.$$
 \end{itemize}
 The second ingredient is the  Todd class of the tangent bundle. This cha\-racteristic class satisfies:
 \begin{itemize}
 \item for a line bundle the Todd class is given by the formula $$td(L)=\frac{c_1(L)}{1-exp(-c_1(L))}=1+\frac{c_1(L)}2+\frac{c_1(L)^2}{12}- \frac{c_1(L)^4}{720} +\frac{c_1(L)^6}{30240}-\dots\,,$$

 \item Todd class is multiplicative $$td(E_1\oplus E_2)=td(E_1) \cup td(E_2)\,.$$
 \end{itemize}

The final elegant form of Riemann--Roch theorem was established by
Alexander Grothendieck, who added the functorial flavour. The Grothen\-dieck-Riemann-Roch deals with
$K(X)$, the $K$-theory of coherent sheaves on $X$. The elements
of $K(X)$ are formal differences of isomorphism classes of coherent sheaves $[E]-[F]$. If
$X$ is smooth, then the {\it coherent sheaves} in the definition
of $K(X)$ can be replaced by the {\it locally free sheaves}, i.e.
holomorphic vector bundles.
 For a proper map of smooth projective varieties $f:X\to Y$ let
  $$f_!(E):=\sum_{i=0}^{\dim(X)}(-1)^iR^if_*(E)\in
K(Y)\,.$$ Here $R^if_*(E)$ is a sheaf with the stalk over $y$
equal to $H^i(f^{-1}(y);E)$.
 In
general $ch(f_!(E))\not=f_*ch(E)$, but   the following diagram
commutes
\begin{theorem}$$\begin{matrix}&&{f_!}\cr
&K(X)&\to&K(Y)\cr
   {ch(-)\cup{ td(TX)}}&\downarrow &&\downarrow&  {ch(-)\cup{ td(TY)}}\cr
&H^*(X)&\to&H^*(Y)\cr
 &&{f_*}\end{matrix}$$
\end{theorem}

The functorial point of view allowed to extend the theory to
singular varieties. Note that in the singular case there is no tangent
bundle $TX$. If the functorial Todd class $td(X)$ existed for
singular $X$ then for an embedding $f:X\hookrightarrow M$ into a
smooth variety  the diagram of Grothendieck-Riemann-Roch would be
commutative:
  $$f_*\left(ch(E)\cup td(X)\right) =ch(f_!(E))\cup td(TM)\,.$$
  For the trivial bundle
$E=\O_X$ whose Chern character is equal to 1, we would obtain.
$$f_*\left(td(X)\right)
 =ch(f_!(\O_X))\cup td(TM)\in H^*(M)\,.$$ In fact this
expression may be taken as the definition of the Todd class, or at
least the definition of its image in $H^*(M)$. But it takes a lot
of work (done by Baum--Fulton--MacPherson, \cite{BFM}) to refine
the construction of the Chern character to obtain the Todd class
localized on $X$, i.e.~to have
 $td^{BFM}(X)\in
H^*(M,M\setminus X)\simeq H_*(X).$ From the practical point of
view the difficulty is hidden in the computation of
 $ch(f_!(\O_X))$.
Not only one has to resolve the sheaf $f_!(\O_X)$ in the sense of
the homological algebra, i.e. find an exact sequence of vector
bundles on $M$:
  $$0\rightarrow E_{\dim(M)}\rightarrow \dots \rightarrow E_1 \rightarrow E_0\rightarrow f_!(\O_X)\rightarrow 0\,.$$
Next one has to apply further nontrivial operations like "graph
construction" to that sequence.

Another, nonequivariant definition of a Todd class for singular
varieties  was proposed by Brasselet-Schuermann-Yokura \cite{BSY}.
Now one imposes a "motivic" condition. One considers not a variety
itself, but a map of varieties. To any map $f:X\to Y$ one
associates a homology class $td(f:X\to Y)\in H_*(Y)$. One demands
that
 \begin{itemize}
\item if $X$ is smooth, then $$td(f:X\to M)=f_*(td(X)\cap[X]\in H_*(M)\,,$$
\item if $Z\subset X$  is a closed subset, $U=X\setminus Z$, then
$$td(f:X\to M)=td(f_{|Z}:Z\to M)+td(f_{|U}:U\to M)\,.$$
\end{itemize}
The existence of  a Todd class is nontrivial, but the computations
are quite effective, provided, that one can resolve the
singularities of $X$ in a geometric sense. We will reserve the
notation $td(-)$ for the motivic Todd class, while the
Baum-Fulton-MacPherson construction will appear only occasionally.

\section{Equivariant version}

We will study the singular subvarieties admitting an action of a
torus. We assume that the torus  $\T=(\C^*)^r$  acts on smooth
complex variety $M$, and $X$ is preserved. For $T=\C^*$ we can say
that we consider quasihomogeneous singularities, but actions of a
bigger  tori appear naturally and then more information about the variety is
involved.  We compute the discussed invariants in the equivariant
cohomology $H^*_\T(M)$ instead of  the usual cohomology. This
cohomology group has reacher structure and it is a refinement of the
usual cohomology. Let us consider cohomology with rational coefficients. For a point we have
$$H^*_\T(pt)=Sym(\T^\vee\otimes\Q)\simeq \Q[t_1,t_2\dots t_r],$$
where $\T^\vee=\Hom(\T,\C^*)$ is the character group. The group
$H^*_\T(M)$ is and algebra  over $H^*_\T(pt)$.

Main advantages  of the equivariant cohomology are
\begin{itemize}  \item If $M$ is a complex smooth compact algebraic variety, then
$$H^*_\T(M)\simeq H^*_\T(pt)\otimes H^*(M)$$
 as modules over $H^*_\T(pt)$, see e.g. \cite{GKM}.

\item Localization theorem:
  the kernel and cokernel of the restriction map  $H^*_\T(M)\to
H^*_\T(M^\T)$ is a torsion
$H^*_\T(pt)$-module, \cite{Qu, GKM}. If $M$ is a complex smooth compact algebraic variety, then the restiction map is in addition injective.
\end{itemize}
The second property is called the Localization Theorem, in fact it
is due to Borel (as pointed by Quillen \cite{Qu}), but it was not formulated by him in that form.
The conclusion is that almost everything about equivariant
cohomology can be read from some   data concentrated at the fixed
points. The tip of this iceberg is the well known fact: the Euler
characteristic $\chi(M)$ is equal to the Euler characteristic of
the fixed point set $\chi(M^\T)$. The situation  when the fixed point set is finite is of particular
interest.

The topological version of the localization theorem was
transformed into a formula which is handy for computations.
Suppose that $M$ is a compact manifold of dimension $n$ and the
torus $(S^1)^r$ acts smoothly with a finite fixed point set. For a
fixed point $p\in M^\T$ define the  Euler class $e(p)\in
H^{2n}_\T(\{p\})$ as the product of characters of $\T$ appearing
in the tangent representation at $p$. The localization theorem
implies that integral of a cohomology class over $M$ can be
expressed by the local data. Precisely

\begin{theorem}[Atiyah-Bott\cite{AB} or Berline-Vergne \cite{BV}]\label{int}
 For a class $a\in H^*_\T(M)$ the integral is the sum of fractions
$$~\int_M a=\sum_{p\in M^\T}
\frac{a_{|p}}{ e(p)}\,.$$
\end{theorem}

\begin{example} \rm Suppose  $M=\P^2$, $T=(\C^*)^3$. Then
 $$M^\T=\{\;p_0=[1:0:0],\; p_1=[0:1:0],\;  p_2=[0:0:1]\;\}$$
 is the fixed point set.
Let $h:=c_1(\O(-1))$ be Chern class of the tautological bundle.
 We apply Berline-Vergne formula to compute the integral for $a=h^2$ $$\int_{\P^2}h^2=\frac{t_0^2}{(t_1-t_0)(t_2-t_0)}+\frac{t_1^2}{(t_0-t_1)(t_2-t_1)}+\frac{t_2^2}{(t_0-t_2)(t_1-t_2)}=$$
$$= -Res_{z=\infty}\frac{z^2}{(z-t_0)(z-t_1)(z-t_2)}$$
If we  set $z=w^{-1}$, then the integral is equal to the
coefficient of  $w$ in $$\frac{{  w}^{-2}}{(w^{-1}-t_0)(w^{-1}-t_1)(w^{-1}-t_2)}=  1$$
as it should be by obvious reasons. See \cite{Zi} for further application of the residue method applied to computations based on Localization Theorem.\end{example}

 We plan to use Localization Theorem to compute some invariants of
$\T$-invariant singular varieties $X\subset M$. The Todd class is
of main interest here. Equally well we  compute Hirzebruch
class $$td_y(X)\in H_*(X)\otimes \Q[y]\,,$$ which is equal to the
formal combination of classes corresponding to sheaves of
differential forms $$\sum_{k=0}^{\dim X}(ch(\Omega^k_X)\cup
td(TX))\cap [X]\,y^k$$ for the smooth $X$ case. The motivic
generalizations of the Hirzebruch class to the singular case was
proven to exist in \cite{BSY}. Also there a relation with other
characteristic classes was discussed. As a particular case $y=-1$
we obtain (after a suitable normalization)
Chern-Schwartz-MacPherson class \cite{McP}.  It is straightforward
to obtain an equivariant version of the Hirzebruch class as it was
done for the Todd class  in \cite{BZ} and  for Chern classes in
\cite{Oh}. These classes live in the {\it equivariant homology}.
Poincar\'e duality isomorphism identifies equivariant homology and
cohomology of $M$. If $M^\T$ has isolated fixed points, then the
class $td_y(X)$ without loss of information can be
replaced by the image in $H^*_\T(M)$. In the notation we indicate
the embedding. We write
 $$td_y(X\to M)\in H^*_\T(M)\,.$$

Now the equivariant Todd class can be studied locally. One can
consider the local Todd class
$$td(X\to M)_{|p}\in \prod_{k=0}^\infty H^k_\T(pt)=\Q[[t_1,t_2,\dots t_r]]$$
for a fixed point $p\in X$.
More generally
$$td_y(X\to M)_{|p}\in \prod_{k=0}^\infty H^k_\T(pt)[y]=\Q[[t_1,t_2,\dots t_r]][y]$$
is an invariant of a $\T$-equivariant singularity germ. Let us
concentrate on the case of $td$ i.e. $y=0$. Due to the integration
formula of Theorem \ref{int} every fixed point of the $\T$-action
contributes to the integral
$$\int_X td(X)=\int_M td(X\to M)=\chi(X;\O_X)\,,$$ that is to the
Todd genus of $X$. The contribution of each point $p\in X^\T$ is
equal to $\frac{td(X\to M )_{|p}}{e(p)}$. Regardless of the
complicated definition the outcome is  just  a Laurent series in
$t_i$. Moreover if we introduce new variables $T_i=e^{-t_i}$ the
answer is of the form
$$
\frac{W(T_1,T_2,\dots T_r)}{\prod_{i=1}^{n}(1-e^{-w_i})}=
\frac{W(T_1,T_2,\dots T_r)}{\prod_{i=1}^{n}(1-\prod_{j=1}^r T_i^{a^j_i})}\,.$$
Here $n$ is the dimension of the ambient space, and $w_i=\sum_{j=1}^ra^j_it_j$
 for $i=1,2,\dots, n$ are the weight
vectors of the tangent representation. The numerator $W(T)=W(T_1,T_2,\dots T_r)$ is a {\it polynomial}
 in $T_i$.
 If the point $p\in X$ is smooth, then
$$W(T_1,T_2,\dots T_r)=1\,,$$
 In general the local
contribution to the $\chi_y$-genus of $X$ is the local Hirzebruch
class
 $$\frac{td_y(X\to M)_{|p}}{e(p)}=\prod \frac {1+y\,e^{-\omega_i}}{1-e^{-\omega_i}}\,,$$
 A question arises for singular points: How to compute the invariant $W(T)$ effectively?
 Using motivic nature of Todd genus one can
decompose $X$ into smooth strata and compute the Todd class
separately for each stratum. Then one has to sum  up the strata
contribution. We will also propose another method, based on
localization theorem, which is effective e.g. for determinant
varieties.

We note that for toric varieties the answer is due to Brion-Vergne
\cite{BV}  and Brylinski-Zhang \cite{BZ}. The Hirzebruch class for
toric varieties is studied in \cite{MS} (although not from the
equivariant point of view) .

\section{Computation by means of a resolution}

Below  we will present some computations just for Todd class, i.e.
for $y=0$.

\begin{example} Whitney umbrella.
\rm We will show step by step how to compute the equivariant Todd class localized at the
origin.\end{example}
Consider the torus $\T=(\C^*)^2$ acting on $\C^3$ by the formula
$$(\,T_1\,,\,T_2\,)\cdot(x_1,x_2,x_3)=(\,T_1T_2x_1\,,\,T_1x_2\,,\,T_2^2x_3\,).$$
 Let us denote the characters $T\to\C^*$
$$ t_1:(T_1,T_2)\to T_1,\quad t_2:(T_1,T_2)\to T_2\,.$$
The considered representation  have the characters
$t_1+t_2\,,\,t_1\,,\,2t_2$.
 The action preserves the Whitney umbrella
 $$X=\{(x_1,x_2,x_3)\in \C^3 \,|\,x_1^2-x_2^2x_3=0\}\,.$$
 We will show that the local equivariant Todd class divided by the Euler class is equal to
$$\frac{td(X\to\C^3)_{|0}}{e(0)}=\frac
       {1+e^{-(t_1+ t_2)}}
  {\left(1-e^{-t_1}\right) \left(1-e^{-2 t_2}\right)}
  $$
We will compute the Todd class  by additivity property.
 Let $$Z=\{(x_1,x_2,x_3)\in \C^3 \,|\, x_1=0,x_2=0\}\quad \text{and}\quad X^o=X\setminus
 Z\,.$$
We have
 $$td(X\to\C^3)=td(X^o\to\C^3)+td(Z\to\C^3)\,.$$
  Let
$$f:\widetilde{X}=\C^2\to \C^3\,,\quad f(u,v)=(uv,u,v^2)\,,$$
be the resolution of the Whitney Umbrella. The map
 $f$  is proper and it is an isomorphism of the open subsets
$$\widetilde{X}^o=(\C\setminus\{0\})\times\C\to X^o.$$
The map $f$ is equivariant provided that the characters of $\C^2$
are $ t_1\,,\,t_2$. Therefore
\begin{align*}td(X\to\C^3)_{|0}&=td(X^o\to\C^3)_{|0}+td(Z\to\C^3)_{|0}\\
&=f_*td(\widetilde{X}^o\to\C^2)_{|0}+td(Z\to\C^3)_{|0}\,.\end{align*}
Let $\C\hookrightarrow\C^2$ be the inclusion as the $v$-coordinate
line. First let us compute
 $$td(\widetilde{X}^o\to\C^2)_{|0}=td(\C^2)_{|0}- td(\C\hookrightarrow\C^2)_{|0}=\frac{
 t_1t_2}{(1-e^{-t_1})(1-e^{-t_2})}-t_1\frac{t_2}{1-e^{-t_2}}\,.$$
The map $f_*:H^k_\T(\C^2)\to H^{k+2}_\T(\C^3)$ is a map of
$H^*_T(pt)$-modules sending $1$ to $\deg(X)=2(t_1+t_2)$, therefore
$td(X\to\C^3)$ is equal to$$2(t_1+t_2)\left(\frac{
 t_1t_2}{(1-e^{-t_1})(1-e^{-t_2})}-t_1\frac{t_2}{1-e^{-t_2}}\right)+(t_1 + t_2)t_1
\frac {2t_2}{1-e^{-2t_2}}=$$

$$=\frac
       {2 t_1t_2(t_1+t_2)\left(1-e^{-2(t_1+ t_2)}\right)}
  {\left(1-e^{-t_1}\right) \left(1-e^{-2 t_2}\right)\left(1-e^{-(t_1+ t_2)}\right)}=td(\C^3)_{|0}\left(1-e^{-2(t_1 t_2)}\right)\,.$$
The image of the Baum-Fulton-MacPherson class is also
equal to $$td(\C^3)_{|0}\cup ch(\O_X)
 =td(\C^3)_{|0}(1-e^{-2(t_1+ t_2)})\,.$$
In general,
if $X$ is a complete intersection $X=\bigcap_{i=1}^k \{f_i=0\}$
then we have
$$ch(\O_X)=\prod_{i=1}^k(1-e^{-\deg(f_i)})\,,$$
where $\deg(f)\in \T^\vee=H^2_\T(pt)$ is understood as the
multidegree of a $\T$-invariant function (if $\T= \C^*$ this is
the quasihomogeneous  degree of $f$). It happens in may cases that
$td(X\to M)=td(M)ch(\O_X)$ but the equality does not hold in
general. An easy counterexample is the cusp $x^3=y^2$ in $\C^2$.
In fact the equality is quite rare. For affine cones over a smooth
hyperplane of degree $d$ in $\P^{n-1}$ the equality holds if and
only if $d\leq n$. One can say that the difference between
Baum-Fulton-MacPherson class and the Brasselet-Schuermann-Yokura
class measures how difficult the singularity is.

There might be a similar relation for full Hirzebruch class. For the
Chern-Schwartz-MacPherson class, i.e. for $y=-1$, the discrepancy between the
expected value and the actual value of the discussed expression
was widely studied. In the case of isolated singularities the
difference is equal up to a sign to the Milnor number, \cite[\S14.1]{Fu}.

\section{Procedure to compute $td(X)$ via localization}

Now, we describe another a method of computing the local
Hirzebruch class. As before we present the computations for Todd
class.
 Suppose that $X\subset M$ is an invariant
subvariety in a $\T$-manifold with isolated fixed points:
$$\int_M td(X\to M)= \sum_{p\in X^\T}\frac{td(X\to M)_{|p}}{e(p)}$$
If $X$  has a decomposition into algebraic
cells the integral is equal to the number of 0--cells.
 Suppose the point $p\in X^\T$ is smooth. Then the contribution
to $\int_Mtd(X\to M)$  is equal to
$$td(X\to M)_{|p}=\prod\frac {{ w_i}}{1-exp(-w_i)}\cdot\prod { n_j}=
{ e(p)}\prod\frac {1}{1-exp(-w_i)}\,,$$
 Here \begin{itemize}
 \item $w_i$ are the
characters  of tangent representation $T_pX$,
\item $n_j$ are the characters  of normal representation.
\end{itemize}
 Suppose that all but one points of $X$ are smooth.
Then one can compute

$$td(X\to M)_{|p_{sing}}=e(p_{sing})\left(\int_M td(X\to M)\;- \sum_{\text{smooth}\; p\in X^\T}\frac{td(X\to M)_{|p}}{e(p)}\right)$$

$$\frac{td(X\to M)_{|p_{sing}}}{e(p_{sing})}=\int_M td(X\to M)\;- \sum_{\text{smooth}\; p\in X^\T}\prod\frac1{1-exp(-w_i)}$$
Amazingly,   computing the local Todd class of a singular variety
this way one does not use its definition but only the knowledge
how it looks like at the smooth points.

 We will show how to  compute the local Todd class when $X$ is the
 set of singular square matrices: $X$ is a hypersurface in
 $M(n\times n;\C)=\C^{n^2}$ given by the equation $\det(A)=0$.
 This set can be compactified: the matrix defines a linear map
 $\C^n\to\C^n$, and its graph is an element of the Grassmannian
 $Gr_n(\C^{2n})$. The closure of the set of singular matrices is
 the Schubert variety of codimension one.

\begin{example}Let $X$ be the codimension one Schubert variety in  $Gr_2(\C^4)$.\end{example}

The torus $T=(\C^*)^4$ acts on $Gr_2(\C^4)$. The action is induced
from the action on $\C^4=\C^2_s\oplus\C^2_t$. The space of linear
maps $\Hom(\C^2_s,\C^2_t)$  is identified with a subset of
$Gr_2(\C^4)$. Let $\overline X$  be the set of planes satisfying the Schubert
condition
 $$\overline X=\{V\in Gr_2(\C^4: V\cap (\C^2_s\oplus 0)\not=0\}\,.$$
 Assume that the action on $\C^2_s$ component is through
 negative characters $-s_1$ and $-s_2$ and on
the second component $\C^2_t$ is through $t_1$ and $t_2$. The
action of $\T$ on $Gr_2(\C^4)$ has exactly six fixed points which are  the coordinate planes.
 The variety $\overline X$ contains all fixed points except
$0\oplus\C_t^2$. The only singular point of $\overline X$ is $p_{sing}=\C^2_s\oplus
0$, which corresponds to $0\in X\subset  \Hom(\C^2_s,\C^2_t)$ in the affine
neighbourhood.  The  action of $\T$ on the tangent space is
through characters $s_i+t_j$. The equation of $X=\overline X\cap
\Hom(\C^2_s,\C^2_t)$ is
$$\det\left(\begin{matrix}a&b\\c&d\end{matrix}\right)=0\,.$$ The
characters of the coordinates are the following:
$$\begin{matrix}\chi_a=s_1+t_1,&\chi_b=s_1+t_2,\\
\chi_c=s_2+t_1,&\chi_d=s_2+t_2.\end{matrix}\,.$$

Let us introduce new variables $S_i=e^{-s_i}$ and $T_i=e^{-t_i}$,
since in the formulas for
$\frac{td(X\to\C^4)_{|p}}{e(p)}=\prod\frac1{1-e^{-x_i}}$ only $e^{-x_i}$
appears, not the character $x_i$ itself.
  We compute the contribution to Berline-Vergne formula at the smooth
  points. For example the fixed point $p=lin\{(1,0,0,0),(0,0,1,0)\}$ contributes with
$$
\frac{1}{\left(1-\frac{S_1}{S_2}\right
   ) (1-\frac1{S_2 T_1})
   \left(1-\frac{T_2}{T_1}\right)}
   = \frac{S_2^2 T_1^2}{(S_1-S_2)
   (1-S_2 T_1) (T_2-T_1)}\,.$$

Now let us compute the Todd class at the singular point:
  $$\frac{td(X\to\C^4)_{|p_{12}}}{e(p_{12})}=\left(\chi(\O_X)
  -\sum_{ \text{smooth}\; p\in X^T}td(X\to\C^4)_{|p}\right)=$$
$$1
  -\frac{S_2^2 T_1^2}{(S_1-S_2)
   (1-S_2 T_1) (T_2-T_1)}
  -\frac{S_1^2 T_1^2}{(S_2-S_1)
   (1-S_1 T_1) (T_2-T_1)}-$$ $$
  -\frac{S_2^2 T_2^2}{(S_1-S_2)
   (1-S_2 T_2) (T_1-T_2)}
  -\frac{S_1^2 T_2^2}{(S_2-S_1)
   (1-S_1 T_2) (T_1-T_2)}=
   $$

    $$= \frac{1-S_1 S_2 T_1
   T_2}{(1-S_1 T_1) (1-S_1
   T_2) (1-S_2 T_1) (1-S_2
   T_2)}$$

Again we see that $$td(X\to\C^4)=td(\C^4)(1-S_1 S_2 T_1
   T_2)=$$ $$=td(\C^4)(1-e^{-(s_1+s_2+t_1+t_2)})=td(\C^4)(1-e^{-\deg(f)})\,,$$
   where $f$ is
   the function defining the variety $X$. The formula for
$td_y(X)_{|p_{12}}$ is much more complicated. On the other hand
the  class $\frac{td_y}{e}$ of the open cell $Gr_2(\C^4)\setminus
X$ is quite friendly
   $$ td(\C^4)(1 + y)^2 \,S_1 S_2 T_1 T_2\Big((1 - y) (1- y\,S_1 S_2 T_1 T_2  ) +
   y (S_1 + S_2) (T_1 + T_2)\Big)$$
Applying this method inductively as described in \cite{We} one can
compute the local Hizebruch class of determinant varieties of
higher dimensions. The result for all $y$ are messy, but for the
special values $y=1$, i.e. the local $L$-class we always have

$$\frac{td_1((\C^{n^2}\setminus X)\to \C^{n^2})}{e(0)}=
2^n\,td(\C^{n^2})\,{\prod_{i<
j}(S_i+S_j)\,\prod_{i<j}(T_i+T_j)}\,{ \prod_i(S_iT_i)}$$ with the
convention $T_i=e^{-t_i}$ and $S_i=e^{-s_i}$. For the  Todd class
in higher dimensions we have
$$\frac{td(X\to\C^{n^2})}{e(0)}=td(\C^{n^2})(1-\prod_{i=1}^n S_iT_i))\,.$$

To finish the list of computations let us display the $td_y$-class
for the open cell for $n=4$ but with the substitution
$s_1=s_2=s_3=s_4=0$, $t_1=t_2=t_3=t_4=t$, that is for the radial
action of $\C^*$ on $\C^4$. The common factor is $ td(\C^{16})(1 +
y)^4T^4$. It should be multiplied by the sum
\begin{align*}
 &\phantom{+T\,}\,\,\,(1-y)^2\, \,(1+y^2)\, \,(1-y+y^2)\,+ \\
&+T\,16 \,(1-y)\, y \,(1-y+y^2)^2\,+ \\
&+T^2\,12 \,(1-y)^2\, y^2 \,(10-13 y+10 y^2)\, + \\
&+T^3\,16 \,(1-y)\, y^2 \,(1+44 y-79 y^2+44 y^3+y^4)\, + \\
&-T^4\,\,(1-y)^2\, y \,(1+29 y+62 y^2-2902 y^3+62 y^4+29 y^5+y^6)\, + \\
&-T^5\,16 \,(1-y)\, y^3 \,(11+62 y-492 y^2+62 y^3+11 y^4)\, + \\
&+T^6\,4 y^3 \,(9-86 y-1139 y^2+3456 y^3-1139 y^4-86 y^5+9 y^6)\, + \\
&+T^7\,16 \,(1-y)\, y^4 \,(11+62 y-492 y^2+62 y^3+11 y^4)\, + \\
&-T^8\,\,(1-y)^2\, y^3 \,(1+29 y+62 y^2-2902 y^3+62 y^4+29 y^5+y^6)\, + \\
&-T^9\,16 \,(1-y)\, y^5 \,(1+44 y-79 y^2+44 y^3+y^4)\, + \\
&+T^{10}\,12 \,(1-y)^2\, y^6 \,(10-13 y+10 y^2)\, + \\
&-T^{11}\,16 \,(1-y)\, y^6 \,(1-y+y^2)^2\, + \\
&+T^{12}\,\,(1-y)^2\, y^6 \,(1+y^2)\, \,(1-y+y^2)\, \end{align*}
 It remains to
say that the meaning of this coefficients is not clear.
When written in that form one can see certain symmetry, which comes from Poincar\'e duality. On the other hand after substitution $y=-1-\delta$, $T=1+S$ one finds that all coefficients are nonnegative. This feature is related to some positivity of the logarithmic vector fields sheaves on the resolution of the determinant variety. In detail it will be explained in \cite{We2}.

Also, by the substitution $T=e^{-\,(y+1)\,t}$ and taking the limit
$y\to -1$ we obtain the Chern-Schwartz-MacPherson class studied in
\cite{AlMi}. The corresponding positivity property in
the nonequivariant case was proven in \cite{Huh} Hopefully computing
the entire Hirzebruch class one may gain better insight into
Chern-Schwartz-MacPherson classes and find nontrivial relation
between them.


\begin{thebibliography}{99}


\bibitem{AlMi} P.~Aluffi,
L.~Constantin Mihalcea, \emph{Chern classes of Schubert cells and
varieties}, J. Algebraic Geom. 18 \,(2009)\,, no. 1, 63–100.

\bibitem{AB} M.~Atiyah, R.~Bott, \emph{The moment map and equivariant
cohomology Topology}, 23 \,(1984)\, 1-28.

\bibitem{BFM}
P.~Baum,  W.~Fulton, R.~MacPherson, \emph{Riemann-Roch for singular
varieties.} Publ. Math., Inst. Hautes Étud. Sci. 45, 101-145
\,(1975)\,

\bibitem{BV} N. Berline, M. Vergne, \emph{Classes caract\'eristiques equivariantes. Formule de
localization en cohomologie \'equivariante}, C.R. Acad. Sc. Paris
295 \,(1982)\,, 539-541.


\bibitem{BZ} J-L. Brylinski,
B. Zhang, \emph{Equivariant Todd Classes for Toric Varieties }
arXiv:math/0311318



\bibitem{BSY} J-P. Brasselet, J. Sch\"urmann, S. Yokura, \emph{Hirzebruch classes and motivic Chern classes for singular spaces.}
J. Topol. Anal. 2, No. 1, 1-55 \,(2010)\,


\bibitem{Fu} W.~Fulton, \emph{Intersection Theory}, Springer 1998

\bibitem{GKM} M.~Goresky, R.~Kottwitz,
R.~MacPherson, \emph{Equivariant Cohomology, Koszul Duality, and
the Localization Theorem}, Invent.~Math.~131, No.1, \,(1998)\,,
25-83

\bibitem{Huh} J.~Huh, \emph{Positivity of Chern classes of Schubert cells and
varieties},  arXiv:1302.5852

\bibitem{McP} R.~MacPherson, \emph{Chern classes for singular algebraic varieties}, Ann. of Math. \,(2)\,, 100, 423-432,
\,(1974)\,

\bibitem{MS} L.~Maxim, J.~Schuerman, \emph{Characteristic classes of singular toric
varieties},  arXiv:1303.4454

\bibitem{Oh} T.~Ohmoto, \emph{Equivariant Chern classes for singular algebraic varieties with
group actions}, Math. Proc. Cambridge Phil. Soc. 140, \,(2006)\,,
115-134


\bibitem{Qu} D.~Quillen, \emph{The Spectrum of an Equivariant Cohomology Ring: I},
Ann. Math., Vol. 94, No. 3, 549-572

\bibitem{Ri} B.~Riemann,  \emph{Theorie der Abel'schen Functionen}. Journal f\"ur die reine und angewandte Mathematik 54 \,(1957)\, 115-155


\bibitem{We} A.~Weber, \emph{Equivariant Chern classes and localization theorem
Journal of Singularities}, Vol. 5 (2012), 153-176

\bibitem{We2} A.~Weber, \emph{Equivariant Hirzebruch class for singular varieties}

\bibitem{Zi} M.~Zielenkiewicz, \emph{Integration over homogenous spaces for classical Lie groups using iterated residues at infinity},	 arXiv:1212.6623, to appear in CEJM

\end{thebibliography}
\end{document}